\documentclass[A4,10pt]{amsart}
\usepackage{amsmath,amssymb,amsrefs,amsfonts,eucal,yfonts,amscd,comment}

\newtheorem{theorem}{Theorem}[section]
\newtheorem{proposition}[theorem]{Proposition}
\newtheorem{lemma}[theorem]{Lemma}
\newtheorem{corollary}[theorem]{Corollary}
\newtheorem{proexample}[theorem]{Example}
\theoremstyle{definition}
\newtheorem{example}[theorem]{Example}
\newtheorem{definition}[theorem]{Definition}
\newtheorem{question}[theorem]{Question}
\newtheorem{problem}[theorem]{Problem}

\renewcommand{\L}{{\mathcal L}}                                 

\newcommand{\Real}{\mathbb R}                                   
\newcommand{\Natural}{{\mathbb N}}                              
\newcommand{\1}{\boldsymbol{1}}

\newcommand{\oto}{\xrightarrow{o}}
\begin{document}
\title{Riesz completions of some spaces of regular operators}

\author{A.W. Wickstead}
\address{Mathematical Sciences Research Centre, Queen's University
Belfast, Belfast BT7 1NN, Northern Ireland.} \email{A.Wickstead@qub.ac.uk}
\begin{abstract}
We  describe the Riesz completion (in the sense of van Haandel) of some spaces of regular operators as explicitly identified subspaces of the regular operators into larger range spaces.
\end{abstract}

 \maketitle
 \section{Introduction}
 The space of regular operators between two Archimedean Riesz spaces is only a Riesz space under the operator order in a few special circumstances, the most well known being when the range space is Dedekind complete (see \cite{Z} for any terms not defined in this paper.) In an effort better to understand partially ordered vector spaces which are not Riesz spaces, van Haandel \cite{vH} introduced the notion of the \emph{Riesz completion} of a  partially ordered vector space. Basically this is the smallest Riesz space into which it may be nicely embedded (precise definitions will be given later.) A very readable account of this and much subsequent related work may be found in \cite{KvG}. In the case of  spaces of regular operators between two Archimedean Riesz spaces, $\L^r(E,F)$, it would be nice to be able to identify the Riesz completion with some other space of operators. One obvious starting point would be to replace the range space $F$ by its Dedekind completion $\hat{F}$ and try to find a realization  of the Riesz completion of $\L^r(E,F)$ as a certain Riesz subspace of $\L^r(E,\hat{F})$. Unfortunately, such a naive approach is doomed to failure as all lattice operations in $\L^r(E,\hat{F})$ satisfy the Riesz-Kantorovich formula but there may be existing lattice operations in $\L^r(E,F)$ which do not, see \cite{E}. It is not difficult to see that this precludes the possibility of realizing the Riesz completion inside $\L^r(E,\hat{F})$ in general. There remains some hope if it is, for some reason, known that any finite suprema that do exist in $\L^r(E,F)$ do satisfy the Riesz-Kantorovich formula.

 In this note we carry out such an identification in a  particularly simple case, namely when the domain is $\ell^\infty_0$, sequences of reals which are eventually constant.  The inspiration for considering this domain came from Example 5.4.15 of \cite{KvG} where there is a description of a vector lattice cover of $\L^r(\ell^\infty_0)$ in terms of matrices representing the operators. Our approach avoids such a representation or any similar concepts and deals solely with operators, in what seems to us a more transparent and natural approach. The same example from \cite{KvG} also prompted us to study the Riesz completion of the space of order continuous operators, in the presence of a mild extra hypothesis on $F$. We conclude the paper by pointing out that many of our results extend to operators from $c$ into a uniformly complete Riesz space and posing some questions concerning dropping the assumption of uniform completeness.

The author would like to thank Michael Elliott,  Onno van Gaans, Anke Kalauch and the referee for correcting several typographical errors and omissions in previous drafts of this paper.

\section{Regular Operators Defined on $\ell^\infty_0$}
Here we gather, for the convenience of the reader, some results from \cite{AW} that we will use. By $\1$ we mean the constantly one sequence and $e_n$ will denote the sequence with all terms zero except for the  $n$'th which is 1. It is clear that $\{\1,e_n:n\in\Natural\}$ is an algebraic basis for $\ell^\infty_0$ as a vector space. $\ell^\infty_0$ is given the pointwise partial order, so that $\{\1,e_n:n\in\Natural\}$ is not a Yudin basis, see \cite{AT}. $F$ will always denote an (at least) Archimedean Riesz space. In spaces $\L(\ell^\infty_0,F)$ we always use the usual  linear operations and pointwise order and $\L^b(\ell^\infty_0,F)$ will denote the order bounded operators in it with $\L^r(\ell^\infty_0,F)$ being the regular operators.

\begin{theorem}[\cite{AG}, Theorem 4.1 and Proposition 4.2]\label{regob}If $F$ is an Archimedean Riesz space then $\L^b(\ell^\infty_0,F)=\L^r(\ell^\infty_0,F)$ which is a Riesz space if and only if $F$ is Dedekind $\sigma$-complete.
\end{theorem}

 The following two propositions are not actually spelt out in \cite{AW}, but are useful in visualising what is happening.
 Recall, first, the following lemma.

\begin{lemma}[\cite{AW}, Lemma 2.2]
If $x_1, x_2, \ldots, x_n$ are arbitrary elements of any Riesz space then
\[\sup\big\{\sum_{k=1}^n\lambda_k
x_k:\lambda_k\in[-1,1]\big\}= \sup\big\{\sum_{k=1}^n\varepsilon_k
x_k:\varepsilon_k\in\{-1,0,1\}\big\}=\sum^n_{k=1} |x_k|\]
\[\mbox{and\ }\sup\big\{\sum_{k=1}^n\lambda_k x_k:\lambda_k\in[0,1]\big\}=
\sup\big\{\sum_{k=1}^n\varepsilon_k
x_k:\varepsilon_k\in\{0,1\}\big\}=\sum^n_{k=1} {x_k}^+.\]
\end{lemma}

\begin{proposition}\label{obexists}Let $F$ be an Archimedean Riesz space and $y_k\in F$ for $k=0,1,2,\dots$. There is $T\in\L^b(\ell^\infty_0,F)$ with $T\1=y_0$ and $Te_n=y_n$ for $n\in\Natural$ if and only if the set $\{\sum_{k=1}^n |y_k|:n\in\Natural\}$ is bounded above in $F$.
\begin{proof}
If such a $T$ exists then the fact that for any $n\in\Natural$ and all choices of $\varepsilon_k\in\{-1,0,1\}$ we have $-\1\le \sum_{k=1}^n\varepsilon_k e_k\le \1$ shows that the set of all possible sums
\[T(\sum_{k=1}^n \varepsilon_k e_k)=\sum_{k=1}^n\varepsilon_k T(e_k)=\sum_{k=1}^n \varepsilon_k y_k\]
is order bounded. By the preceding lemma, the set $\{\sum_{k=1}^n |y_k|:n\in\Natural\}$ is order bounded.

Now suppose that $\{\sum_{k=1}^n |y_k|:n\in\Natural\}$ is bounded above in $F$ by $u$. It is now routine to define a linear operator $T:\ell^\infty_0\to F$ by $T(x_k)=T(x_n\1+\sum_{k=1}^n (x_k-x_n)e_k=x_n y_0+\sum_{k=1}^n(x_k-x_n)y_k$,where $x_k=x_n$ for $k\ge n$.

In order to show that $T$ is order bounded, it suffices to show that the image under $T$ of the order  interval $[-\1,\1]$ in $\ell^\infty_0$ is order bounded in $F$. If $x\in [-\1,\1]$ we may write
$x=x_n\1+\sum_{k=1}^n (x_k-x_n)e_k$ where $x_k\in [-1,1]$ for $1\le k\le n$, so that $x_k-x_n\in [-2,2]$, for $1\le k\le n$. Now
\[|Tx|=|x_nT\1+\sum_{k=1}^n (x_k-x_n) Te_k|=|x_n y_0+\sum_{k=1}^n (x_k-x_n) y_k|\le |y_0|+\sum_{k=1}^n 2|y_k|\le |y_0|+2u,\]
so that $T([-\1,\1]\subseteq[-|y_0|-2u,|y_0|+2u]$ and $T$ is certainly order bounded.
\end{proof}
\end{proposition}

 \begin{proposition}\label{oborder}Let $F$ be an Archimedean Riesz  space and $T\in \L^b(\ell^\infty_0,F)$. $T$ is positive if and only if (i) $Te_n\ge 0$ for all $n\in\Natural$ and (ii) $T\1\ge \sum_{k=1}^n Te_k$ for all $n\in\Natural$.
\begin{proof}The proof of ``only if'' is clear as each $e_k\ge 0$ and $\1\ge \sum_{k=1}^n e_k$ for all $n\in\Natural$.

Suppose that inequalities (i) and (ii) are satisfied. Let $x\in \ell^\infty_{0+}$, the positive cone of $\ell^\infty_0$. There is $n\in\Natural$ such that
\begin{align*}x=(x_1,x_2,\dots,x_n,x_n,\dots)&=x_n\1+\sum_{k=1}^{n-1}(x_k-x_n)e_k\\
&= x_n[\1-\sum_{k=1}^{n-1} e_k]+\sum_{k=1}^{n-1} x_ke_k.
\end{align*}
Apply $T$, bear in mind (i) and (ii) plus the fact that each $x_k\ge 0$ and we see that
\[Tx=x_n T(\1-\sum_{k=1}^n e_k)+\sum_{k=1}^{n-1} x_k e_k\ge 0.\]
\end{proof}
\end{proposition}
 \begin{theorem}[\cite{AW}, Theorems 4.4 and 4.6]\label{pospart}If $F$ is an Archimedean Riesz space and $T\in\L^r(\ell^\infty_0,F)$ then $T^+$ exists if and only if the supremum
 \[\sup\{[\sum_{k=1}^n (Te_k)^+]\vee[T\1+\sum_{k=1}^n (Te_k)^-]:n\in\Natural\}\]
 exists in $F$. Furthermore, in such cases, $T^+$ satisfies the Riesz-Kantorovich formula.
 \end{theorem}

\section{The Riesz completion of $\L^r(\ell^\infty_0,F)$}
In \cite{KvG}, Theorem 2.4.5, one may find a proof that every pre-Riesz space has a Riesz completion. For our purposes the definition of a pre-Riesz space is not important. All we need to know is that every directed, i.e. positively generated, Archimedean partially ordered vector space is pre-Riesz, \cite{KvG} Proposition 2.2.3. All our spaces of operators are pre-Riesz. We do need to know the definition of the Riesz completion.

Recall that a linear subspace $E$ of a partially ordered vector space $F$ is \emph{order dense} if for every $y\in F$ we have $y=\inf\{x\in E:x\ge y\}$ or, equivalently, $y=\sup\{x\in E:x\le y\}$.

\begin{definition}If $E$ is a pre-Riesz space then a Riesz space $F$ is a \emph{vector lattice cover} of $E$ if there is a bipositive linear map $i:E\to F$ such that $i(E)$ is order dense in $F$. If also $i(E)$ generates $F$ as a Riesz space then we call $F$ the \emph{Riesz completion} of $E$.
\end{definition}

The Riesz completion is essentially unique and we may clearly think of the Riesz completion $F$ as containing $E$. We follow \cite{KvG} by denoting the Riesz completion of $E$ by $E^\rho$.

\begin{definition}Let $E$ and $F$ be Archimedean Riesz spaces. We say that $E$ is \emph{sequentially order closed} in $F$ if
\begin{enumerate}
\item $E$ is an order dense vector sublattice of $F$.
\item Any increasing sequence in $E$ with a supremum in $E$ has the same supremum in $F$.
\item If $(x_n)$ is an increasing sequence in $E$ which is bounded above then $(x_n)$ has a supremum in $F$.
\end{enumerate}
We say that $F$ is a \emph{sequential order closure} of $E$ if also
\begin{enumerate}
\item[(4)] Condition (3) does not hold for any smaller vector sublattice of $F$.
\end{enumerate}
\end{definition}

Note that there is no implication that $F$ is Dedekind $\sigma$-complete. Indeed many of our results would become almost trivial if that were possible.

\begin{proposition}Every Archimedean vector lattice $E$ has a sequential order closure.
\begin{proof}
Take $F$ to be the Dedekind completion of $E$. Let
$\overline{F}$ be the set of all suprema (in $F$) of increasing order bounded sequences in $E$ and define $G=\overline{F}-\overline{F}$.
Once we show that $G$ is a vector lattice, conditions (1), (2), (3) and (4) will be obvious.

It should be clear that $\overline{F}$ is closed under addition and multiplication by non-negative reals.\footnote{Such sets are known by various names. The term \emph{wedge} has the advantage of non-ambiguity. Some authors would say that it is a cone, although most authors would insist that a cone $C$ satisfy $C\cap (-C)=\{0\}$. In our case it can even happen that $\overline{Y}=-\overline{Y}$, for example if $X$ is almost Dedekind $\sigma$-complete.} It is then routine to show that $G=\overline{F}-\overline{F}$ is a vector subspace of $F$. In order to show that $G$ is a sublattice of $F$, it suffices to show that if $x,y\in \overline{F}$ then $(x-y)^+\in G$. It is trivial to check that $x\vee y\in \overline{F}$ and then we need only use the identity $(x-y)^+=x\vee y-y$.
\end{proof}
\end{proposition}

It is clear that the sequential order closure of an Archimedean vector lattice is essentially  unique.

We write $E^\sigma$ for the (essentially unique) sequential order closure of $E$.

Recall from \cite{AL1} that a Riesz space $E$ is \emph{almost Dedekind $\sigma$-complete}\footnote{They actually use the phrase \emph{almost $\sigma$-Dedekind complete} but we feel that our terminology is rather more standard.} if it can be embedded as a super order dense Riesz subspace of a Dedekind $\sigma$-complete Riesz space $F$, i.e. every $y\in F$ is both the supremum and the infimum of suitably selected \emph{sequences} from  $E$. Such an $F$ is termed the Dedekind $\sigma$-completion of $E$. In particular if $E$ is order separable then it is almost Dedekind $sigma$-complete.  It is evident that, in this case, $E^\sigma=F$ and therefore $E^\sigma$ is Dedekind $\sigma$-complete and hence is uniformly complete, which will be needed later. In this case $E^\sigma$ coincides with the  Dedekind $\sigma$-completion\footnote{The reader is cautioned that van Haandel's Dedekind $\sigma$-completion in \cite{vH} is a different object.} of $E$ in the sense of Quinn \cite{Q} and of Aliprantis and Langford \cite{AL2}. Example (iv) on pages 424--425 of \cite{AL1} is an Archimedean Riesz space whch has the principal projection property but is not almost Dedekind $\sigma$-complete. This example is, of course, far from being uniformly complete. Given the importance of this notion in our work, we include an example, due to Eugene Bilokopytov and included here with his permission, of a Banach lattice that is not almost Dedekind $\sigma$-complete. Michael Elliott has kindly pointed out to the author that Theorem 4.1 of \cite{AW2} describes two families of unital AM-spaces which fail to be almost Dedekind $\sigma$-complete.

\begin{example} Let  $I$ be an uncountable set endowed with the discrete topology ad let $\alpha I=I\cup\{\infty\}$ be its one point compactifaction. A real-valued function $f$ on $\alpha I$ lies in $C(\alpha I)$ precisely when $\forall\epsilon>0$ the set $\{i\in I:|f(i)-f(\infty)|<\epsilon\}$ is finite, so in particular
\[\{i\in I:f(i)\ne f(\infty)\}=\bigcup_{n=1}^\infty\{i\in I:|f(i)-f(\infty)|>\frac1n\}\]
is countable. Suppose that $C(\alpha I)$ were almost Dedekind $\sigma$-complete then we may embed it as a super order dense sublattice of a Dedekind $\sigma$-complete Riesz space $F$. We use $\chi_A$ to denote the characteristic function of a set $A\subseteq I$, which need not lie in $C(\alpha I)$, and $\chi_a$ for the characteristic fnction of the set $\{a\}$ for $a\in I$, which wil lie in $C(\alpha I)$. Pick a countably infinite set $A\subset I$ and consider the family $\{\chi_a:a\in A\}$ in $C(\alpha I)$. This is bounded above by the constantly one function, so has a supremum $y\in F$. If $b\in I\setminus A$ then $\chi_a\perp \chi_b$ for all $a\in A$, so that $y\perp \chi_b$.

By the super order density, $y$ is the infimum in $F$ of some sequence $f_n\in C(\alpha I)$. As each $f_n\ge y\ge \chi_a$ for all $a\in A$, we have $f_n\ge \chi_A$ and therefore $f_n(\infty)\ge 1$. Hence, for each $n\in \Natural$, there s a finite set $F_n$ such that $f_n(i)\ge \frac12$ for all $i\in I\setminus F_n$. Thus for all $n\in \Natural$ we have $f_n(i)\ge\frac12$ for all $i\in I\setminus \bigcup_{n=1}^\infty F_n$. Pick any $b\in I\setminus (A\cup\bigcup_{n=1}^\infty F_n)$, which is certainly non-empty, and we have $f_n\ge \frac12\chi_b$ for all $n\in\Natural$ and hence $y=\bigwedge_{n=1}^\infty f_n\ge \frac12 \chi_b$, which contradicts $y\perp \chi_b$.

\end{example}

\begin{theorem}\label{Rcompl}
If $F$ is an Archimedean Riesz space then the Riesz completion of $\L^r(\ell^\infty_0,F)$, $\L^r(\ell^\infty_0,F)^\rho$, may be identified with the space of operators $\{T\in \L^r(\ell^\infty_0,F^\sigma):Te_n\in F\ \forall n\in\Natural\}$.
\begin{proof}
In this proof we will denote $\{T\in\L^r(\ell^\infty_0,F^\sigma): Te_n\in F\ \forall n\in\Natural\}$ by $H$, for brevity.

Our first task is to show that $H$ is a Riesz space. By Proposition \ref{obexists} applied in $\L^r(\ell^\infty_0,F^\sigma)$, an operator in $H$ is specified by $T\1\in F^\sigma$ and $Te_n\in F$ provided only that $\{\sum_{k=1}^n |Te_k|:n\in\Natural\}$ is bounded above in $F^\sigma$ and therefore in $F$. By Theorem \ref{pospart} this has a positive part in $\L^r(\ell^\infty_0, F^\sigma)$ provided that
\[\sup\{[\sum_{k=1}^n (Te_k)^+]\vee [T\1+\sum_{k=1}^n (Te_k)^-:n\in\Natural\}\]
exists in $F^\sigma$. But this supremum is equal to
\begin{align*}
\sup\{\sum_{k=1}^n (Te_k)^+:n\in\Natural\}&\vee \sup\{\{T\1+\sum_{k=1}^n (Te_k)^-:n\in\Natural\}\\
=\sup\{\sum_{k=1}^n (Te_k)^+:n\in\Natural\}&\vee [T\1+\sup\{\sum_{k=1}^n (Te_k)^-:n\in\Natural\}],
\end{align*}
and $\sup\{\sum_{k=1}^n (Te_k)^\pm:n\in\Natural\}$ are countable order bounded sets from $F$ so have suprema in $F^\sigma$ whilst certainly $T\1\in F^\sigma$. Thus $T^+$ exists in $\L^r(\ell^\infty_0,F^\sigma)$ for all $T\in H$. As, for each $n\in\Natural$, $T^+e_n=(Te_n)^+\in F$, we see that $T^+\in H$  which is enough to show that $H$ is a Riesz space.

The next step is to show that $\L^r(\ell^\infty_0,F)$ is order dense in $H$. If $T\in H$, then $T\1=\sup\{y\in F:y\le T\1\}$, as $F$ is order dense in $F^\sigma$. For each $y\in F$ with $y\le T\1$, define $T_y\in\L^r(\ell^\infty_0,F)$ with $T_ye_n=Te_n$ for $n\in\Natural$ and $T_y\1=y$. We have $T_y\le T$ and $T$ is the supremum in $H$ of the $T_y$ as if $T_y\le S\le T$ for all these $y$, then $T_ye_n=Te_n\le  Se_n\le Te_n$ for all $n\in\Natural$ and $T\1\ge S\1\ge \sup\{T_y\1:y\le T\1\}=\sup\{y:y\le T\1\}=T\1$ so that $S=T$.

Finally, we need to show that $\L^r(\ell^\infty_0,F)$ generates $H$ as a Riesz space. If $T\in H$ then there is $T_0\in \L^r(\ell^\infty_0,F)$ with $T_0 e_n=Te_n$ for all $n\in\Natural$ and $T_0\1=0$, as the order boundedness criterion is the same in $\L^r(\ell^\infty_0,F)$ as in $H$. Thus it suffices to show that operators $T\in H$ for which $Te_n=0$ for all $n\in\Natural$ lie in the sublattice of $H$ generated by $\L^r(\ell^\infty_0,F)$. Such operators are completely specified by the value of $T\1\in F^\sigma$. A typical element of $F^\sigma$ is $t=u-v$ where both $u$ and $v$ are a countable supremum from $F$. Suppose that  $u=\bigvee_{n=1}^\infty y_n$ from $F$, which we may assume is increasing. Let $x_n=y_n-y_1$, so that each $x_1=0$. Define $U_1\in \L^r(\ell^\infty_0,F)$ with $U_1\1=0$ and $U_1e_n=x_n-x_{n+1}$. By Proposition \ref{pospart} the positive part of  $U_1$ in $H$ satisfies $U_1^+e_n=(U_1e_n)^+=(x_n-x_{n+1})^+=0$, whilst
\begin{align*}
U_1\1&=\sup\{[\sum_{k=1}^n (U_1e_k)^+]\vee[U_1\1+\sum_{k=1}^n (U_1e_k)^-]:n\in\Natural\}\\
&=\sup\{[\sum_{k=1}^n (x_k-x_{k+1})^+]\vee[0+\sum_{k=1}^n  (x_k-x_{k+1})^-]:n\in\Natural\}\\
&=\sup\{\sum_{k=1}^n x_{k+1}-x_k:n\in\Natural\}=\sup\{x_{n+1}:n\in\Natural\}=u-y_1.
\end{align*}
If $U_2\in \L^r(\ell^\infty_0,F)$ with $U_2e_n=0$ for all $n\in\Natural$ and $U_2\1=y_1\in F$ then we see that, setting $U=U_1^++U_2$, we have $U\1=u$ and $Ue_n=0$ for $n\in\Natural$. Similarly we may find $V_1,V_2\in \L^r(\ell^\infty_0,F)$ such that $V=V_1^++V_2$ has $V\1=v$ and $Ve_n=0$ for all $n\in\Natural$. Now $T=U-V=U_1^++U_2-V_1^+-V_2$ is certainly in the Riesz subspace of $H$ generated by $\L^r(\ell^\infty_0,F)$ and  has $T\1=U\1-V\1=u-v=t$ and $Te_n=0$ for  all $n\in\Natural$, completing the proof.
\end{proof}
\end{theorem}

A strong, and very useful property  that a pre-Riesz space may have is that of being \emph{pervasive}.

\begin{definition} A pre-Riesz space $E$ is \emph{pervasive} if whenever $0\ne y\in (E^\rho)_+$ there is $0\ne x\in E_+$ with $x\le y$.
\end{definition}

\begin{theorem}\label{pervasive}If $F$ is an Archimedean Riesz space then $\L^r(\ell^\infty_0,F)$ is pervasive.
\begin{proof}
We start with a non-zero positive operator $T\in \L^r(\ell^\infty_0,F)^\rho$ which we identify with $\{T\in \L^r(\ell^\infty_0,F^\sigma):Te_n\in F\ \forall n\in\Natural\}$. We consider two cases. Suppose first that $Te_n=0$ for all $n\in\Natural$. Certainly $T\1\ne 0$ and as $F$ is order dense in $F^\sigma$ there is $0\ne y\in F_+$ with $y\le T\1$. By Proposition \ref{obexists} there is $S\in\L^r(\ell^\infty_0,F)$ with $S\1=y$ and $Se_n=0$ for all $n\in\Natural$. It is simple to see that $0\le S\le T$ and that $S\ne 0$. In the second case, one of the $Te_n\ne 0$. Without loss of generality, let us assume that $Te_1\ne 0$. Using Proposition \ref{obexists} there is $S\in \L^r(\ell^\infty_0,F)$ with $S\1=Se_1=Te_1$ and $Se_n=0$ for $n>1$. Certainly $0\le S\ne 0$, and $(T-S)e_n\ge 0$ for all $n\in \Natural$. Also, for each $n\in\Natural$
\[(T-S)\1=T\1-Te_1\ge \sum_{k=2}^n Te_k=\sum_{k=1}^n (T-S)e_k,\]
so that $S\le T$ as required.
\end{proof}
\end{theorem}

\section{Order continuous operators defined on $\ell^\infty_0$}

A net $(x_\gamma)_{\gamma\in \Gamma}$ in an Archimedean Riesz space $E$ is \emph{order convergent} to $x\in E$ if there is a net $(a_\lambda)_{\lambda\in \Lambda}$ in $E$ such that $a_\lambda\downarrow 0$ and for all $\lambda\in\Lambda$ there is $\gamma_0\in\Gamma$ such that, for all $\gamma\ge\gamma_0$, $|x_\gamma-x|\le a_\lambda$. In this case we write $x_\gamma\oto x$. Similarly, but rather more simply, a sequence $(x_n)$ in $E$ is \emph{order convergent} to $x\in E$ if there is a sequence $(a_n)$ in $E$ such that $a_n\downarrow 0$ and $|x_n-x|\le a_n$ for all $n\in\Natural$. Again, we write $x_n\oto x$.

If $E$ and $F$ are Archimedean Riesz spaces then a linear operator $T:E\to F$ is \emph{order continuous} if whenever $(x_\alpha)_{\alpha\in A}$ is a net in $E$ such that $x_\alpha\oto 0$ then $Tx_\alpha\oto 0\in F$. Similarly, $T:E\to F$ is \emph{$\sigma$-order continuous} if this holds for sequences. For positive operators, there is a simpler characterization that $T$ is order continuous if and only if $Tx_\alpha\downarrow 0$ whenever $x_\alpha\downarrow  0$, with a similar result for $\sigma$-order continuity.

By Theorem 2.1 and Remark 2.2 of \cite{AS}, every order continuous operator is order bounded, but $\sigma$-order continuous operators need not be order bounded. Accordingly, we write $\L^n(E,F)$ (resp. $\L^c(E,F)$) for the set of order continuous (resp. $\sigma$-order continuous and  order bounded) operators from $E$ into $F$. The fundamental properties of order convergence in Riesz spaces make it clear that both $\L^n(E,F)$ and $\L^c(E,F)$ are linear subspaces of $\L^b(E,F)$. At present we are unable to say much more in general. This section is devoted to investigating order continuous operators defined on $\ell^\infty_0$, before attempting, with only partial success, to obtain a  description of the Riesz completion. Note that in view of Theorem \ref{regob} order continuous operators defined on $\ell^\infty_0$, being order bounded, are automatically regular. This does not, however, mean that every order continuous operator is the difference of two positive order continuous operators.

We make a start on describing order continuous operators on $\ell^\infty_0$ with a simple necessary condition.

\begin{proposition}\label{posoc} Let $F$ be an Archimedean Riesz space and $T\in \L^r(\ell^\infty_0,F)_+$ then  the following are equivalent:
\begin{enumerate}
\item $T$ is order continuous.
\item $T$ is $\sigma$-order continuous.
\item $T\1-\sum_{k=1}^n Te_k\downarrow 0\in F$.
\end{enumerate}
\begin{proof}
Certainly (1) implies (2). Assuming (2), take $x_n=\1-\sum_{k=1}^n e_k$ and observe that $x_n\downarrow 0$ in $\ell^\infty_0$ so $Tx_n\downarrow 0$ in $F$, giving us (3) Now assume (3). Let $(x^\alpha)_{\alpha\in A}$ be a positive net in $\ell^\infty_0$ decreasing to 0 and set $x^\alpha=(x^\alpha_1, x^\alpha_2,\dots)$. From Proposition \ref{oborder} the fact that $x^\alpha\downarrow 0\in c$ implies that for each $n\in\Natural$ we have $x^\alpha_n\downarrow 0\in\Real$ and also that $(\lim_{n\to\infty} x^\alpha_n)_{\alpha\in A}\downarrow$ (but not necessarily to 0.) Pick any $\alpha_0\in A$ and let $M=\|x^{\alpha_0}\|_\infty$. For any $\epsilon>0$ and $n\in\Natural$, we may find $\alpha_1,\alpha_2,\dots,\alpha_n\ge \alpha_0$ such that $0\le x^{\alpha_k}_k<\epsilon$. Now choose $\beta_{\epsilon,n} \ge \alpha_0,\alpha_1,\dots\alpha_n$, so that $0\le x^{\beta_{\epsilon,n}}_k<\epsilon$ for $1\le k\le n$ and $0\le x^{\beta_{\epsilon,n}}_m\le x^{\alpha_0}_m\le M$ for all $m\in\Natural$. It follows that
\[0\le x^{\beta_{\epsilon,n}} \le \epsilon\1+(0,0,\dots,0,M,M,\dots)=\epsilon\1+M(\1-\sum_{k=1}^n e_k).\]
Hence $0\le T(x^{\beta_{\epsilon,n}})\le \epsilon T\1+M(T\1-\sum_{k=1}^n Te_k)$ so that
\[0\le \bigwedge_{\alpha\in A} T(x^\alpha)\le \bigwedge_{n\in \Natural,\epsilon>0} \epsilon T\1+M(T\1-\sum_{k=1}^n Te_k)=0.\]
I.e. $Tx^\alpha\oto 0$. The extension to an arbitrary net $(x^\alpha)_{\alpha\in A}$ for which $x_\alpha\oto 0$ is routine.
\end{proof}
\end{proposition}

We deal next with the special case of operators that have a positive part, which is fairly easy to obtain. We need two simple preliminary results.

\begin{lemma}\label{uppersets}Let $F$ be an Archimedean Riesz space, and $(x_n), (y_n)$ be increasing sequences in $F$ with $x_n-y_n\oto 0$. The three sets $A=\{z\in F:z\ge x_n\ \forall n\in\Natural\}$, $B=\{z\in F:z\ge y_n\ \forall n\in\Natural\}$ and $C=\{z\in F:z\ge x_n\vee y_n\ \forall n\in\Natural\}$ are all equal.
\begin{proof}
Clearly $C\subseteq A,B$. We prove that $A\subseteq C$, with the proof that $B\subseteq C$ being virtually identical. Suppose that $z\in A\setminus C$, then $x_n\le z$ for all $n\in\Natural$, but there is $m\in\Natural$ with $x_m\vee y_m\not\le z$, so that $y_m\not\le z$, so that $0\le (y_m-z)^+\ne 0$. If $n\ge m$  then
\[(y_m-z)^+\le (y_n-z)^+\le (y_n-x_n)^+\le |y_n-x_n|,\]
which contradicts $x_n-y_n\oto 0$.
\end{proof}
\end{lemma}

\begin{corollary}\label{upperlimits}Let $F$ be an Archimedean Riesz space, and $(x_n), (y_n)$ be increasing sequences in $F$ with $x_n\vee y_n\uparrow \ell\in F$. If $x_n-y_n\oto 0$ then $x_n\uparrow \ell$ and $y_n\uparrow \ell$.
\begin{proof}
With the notation of the preceding lemma, $C$ has a least element, $\ell$, hence so have $A$ and $B$, from which the claim is immediate.
\end{proof}
\end{corollary}

\begin{theorem}\label{ocifpospart}
Let $F$ be an Archimedean Riesz space and let $T\in\L^r(\ell^\infty_0,F)$ have a positive part, then the following are equivalent:
\begin{enumerate}
\item $T$ is order continuous.
\item $T$ is $\sigma$-order continuous.
\item $T^+$ and $T^-$ are order continuous.
\item $T\1-\sum_{k=1}^n Te_k\oto 0$.
\end{enumerate}
\begin{proof}
We already know that (1) implies (2), which as in the proof of Proposition \ref{posoc} implies (4), and also that (4) implies (1) if we also assume that $T$ is positive. It is also clear that (3) implies (2). The collection of operators $T$ in $\L^r(\ell^\infty_0,F)$ which satisfy (4) is clearly a linear subspace of $\L^r(\ell^\infty_0,F)$, which we will temporarily denote by $\L^\diamond(\ell^\infty_0,F)$. As $T$ has a positive part, Theorem 4.4 of \cite{AW} tells us that $T^+e_n=(Te_n)^+$ for $n\in\Natural$, and that
\begin{align*}
T^+\1&=\sup\{[\sum_{k=1}^n (Te_k)^+]\vee[T\1+\sum_{k=1}^n (Te_k)^-]:n\in\Natural\}
\\&=\sup\{\sum_{k=1}^n (Te_k)^+:n\in\Natural\}\vee \sup\{T\1+\sum_{k=1}^n (Te_k)^-:n\in\Natural\}
\end{align*}
with this supremum existing in $F$. By (4),
\[[T\1+\sum_{k=1}^n (Te_k)^-]-[\sum_{k=1}^n (Te_k)^+]=T\1-\sum_{k=1}^n Te_k\oto0,\eqno{(*)}\]
so Corollary \ref{upperlimits} tells us that $\sum_{k=1}^n (Te_k)^+=\sum_{k=1}^n T^+e_k\uparrow T^+\1.$ Hence
 $T^+\in \L^\diamond(c,F)$, and certainly $0\le T^-=(T^+-T)\in \L^\diamond(\ell^\infty_0,F)$. By  Proposition \ref{posoc}, $T^+$ and $T^-$ are order continuous, establishing (3).
\end{proof}
\end{theorem}

\begin{corollary}\label{ocifDsc}
Let $F$ be a  Dedekind $\sigma$-complete Riesz space then $\L^r(\ell^\infty_0,F)$ is a Riesz space and for $T\in\L^r(\ell^\infty_0,F)$, the following are equivalent:
\begin{enumerate}
\item $T$ is order continuous.
\item $T$ is $\sigma$-order continuous.
\item $T\1-\sum_{k=1}^n Te_k\oto 0$.
\end{enumerate}
\end{corollary}

It is simple to see that if $F$ is Dedekind $\sigma$-complete then $\L^n(\ell^\infty_0,F)$ is a vector lattice ideal in $\L^r(\ell^\infty_0,F)$. A little extra work shows that it is actually a band, however direct construction of the associated band projection allows us to show that it is actually a projection band, replicating the conclusion of Ogasawara's theorem for operators into a Dedekind complete Riesz space.

\begin{theorem}\label{ocband}If $F$ is a Dedekind $\sigma$-complete Riesz space then
\begin{enumerate}
\item $\L^n(\ell^\infty_0,F)$ is a projection band in $\L^r(\ell^\infty_0,F)$.
\item The band projection $P$ onto $\L^n(\ell^\infty_0,F)$ maps the operator $T$ to an operator $P(T)$ where $P(T)e_n=Te_n$ for all $n\in\Natural$ and $\sum_{k=1}^n Te_n\oto P(T)\1$.
\item The complementary band, $\L^n(\ell^\infty_0,F)^d$, consists of all operators of the form $x=(x_n)\mapsto \lim_{n\to\infty}x_n\times y$ for $y\in F$, and hence
\item $\L^n(\ell^\infty_0,F)^d$ is linearly order isomorphic to $F$.
\end{enumerate}
\begin{proof}
If $T\in\L^r(c,F)_+$ then $\sum_{k=1}^n Te_k$ is bounded above, as $n$ ranges over the natural numbers, by $T\1$ and has a supremum $u$ as $F$ is Dedekind $\sigma$-complete. Taking $P(T)e_n=Te_n$ for $n\in\Natural$ and $P(T)\1=u$ then extending $P(T)$ linearly to lie in $\L^r(\ell^\infty_0,F)$ we see immediately that $P(T)\in \L^n(\ell^\infty_0,F)$. $P$ is additive and positively homogeneous on $\L^r(\ell^\infty_0,F)_+$ so extends to a linear operator $P:\L^r(\ell^\infty_0,F)\to \L^n(\ell^\infty_0,F)$. It is clear that $P$ is positive, that $P(T)\le T$ if $T$ is positive and that $P^2=P$. To complete the proof of (1) and (2) we need only note that if $\sum_{k=1}^n T_1e_k\uparrow s$  and $\sum_{k=1}^n T_2e_k\uparrow t$ then $\sum_{k=1}^n (T_1-T_2)e_k\oto s-t$.

A typical element of $\L^n(\ell^\infty_0,F)^d$ is $S=T-P(T)$ for $T\in\L^r(\ell^\infty_0,F)$. But $Se_n=Te_n-P(T)e_n=0$ for all $n\in\Natural$ and $S\1=T\1-\sum_{k=1}^\infty Te_k\in F$, where the series is order convergent. This establishes (3) and then (4) follows immediately.
\end{proof}

\end{theorem}

\begin{theorem}\label{ocasDc}
Let $F$ be an almost Dedekind $\sigma$-complete Riesz space and $T\in \L^r(\ell^\infty_0,F)$, then the following are equivalent:
\begin{enumerate}
\item $T$ is order continuous.
\item $T$ is $\sigma$-order continuous.
\item $T\1-\sum_{k=1}^n Te_k\oto 0$.
\end{enumerate}
\begin{proof}
As in the proof of Theorem \ref{ocifpospart}, we need only prove that (3) implies (1), and we retain the notation of that proof. By $(*)$ and Lemma \ref{uppersets}, the two sets $A=\{\sum_{k=1}^n (Te_k)^+:n\in\Natural\}$ and $B=\{T\1+\sum_{k=1}^n (Te_k)^-:n\in\Natural\}$ have the same upper bounds. Let $u$ be one such upper bound. As $F$ is almost Dedekind $\sigma$-complete, we may embed it in a Dedekind $\sigma$-complete Riesz space $F^\sigma$ with every $y\in F^\sigma$ being both the infimum and supremum of a sequence from $F$. Let $x$ be the supremum in $F^\sigma$ of the set $A$ and consider $u-x\in F^\sigma_+$. This is the supremum of a sequence $(y_n)_{n=1}^\infty$ in $F$, which we may assume is increasing and positive. Define also $y_0=0$ and now set $z_n=y_n-y_{n-1}$ for $n\in\Natural$. We now define, using Proposition~\ref{obexists}, $S\in \L^r(\ell^\infty_0,F)$ by $S\1=u$ and $Se_n=(Te_n)^++z_n$ for $n\in\Natural$. This exists as
\begin{align*}
\sup\{\sum_{k=1}^n ((Te_k)^++z_k):n\in\Natural\}&=\sup\{y_n+\sum_{k=1}^n(Te_k)^+:n\in\Natural\}\\
&=\sup\{y_n:n\in\Natural\}+\sup\{\sum_{k=1}^n (Te_k)^+:n\in\Natural\}\\
&=(u-x)+x=u
\end{align*}
in $F$ and hence in $F^\sigma$. Clearly, $S$ is positive and this equality states precisely that $\sum_{k=1}^n Se_k\uparrow S\1$ so that $S\in\L^\diamond(\ell^\infty_0,F)_+$ and therefore $S\in\L^n(\ell^\infty_0,F)$. We also have $S\ge T$ as $Se_n=(Te_n)^++z_n\ge Te_n$ for all $n\in\Natural$ and, for any $n\in\Natural$,
\begin{align*}
S\1-\sum_{k=1}^n Se_k&=u-[y_n+\sum_{k=1}^n (Te_k)^+]\\
&\ge [y_n+T\1+\sum_{k=1}^n (Te_k)^-]-[y_n+\sum_{k=1}^n(Te_k)^+]\\
&=T\1-\sum_{k=1}^n Te_k.
\end{align*}
 Hence $(S-T)\in \L^\diamond(\ell^\infty_0,F)_+$. By  Proposition \ref{posoc}, we  again see that $S-T$ is order continuous and hence so is $T$.
\end{proof}
\end{theorem}

Recall that a partially ordered vector space $E$ is \emph{positively generated} if every $x\in E$ may be written as $x=y-z$ with $y,z\in E_+$. The last proof also establishes:

\begin{theorem}\label{ocpg}
If $F$ is an almost Dedekind $\sigma$-complete Riesz space then $\L^n(\ell^\infty_0,F)$ is positively generated.
\end{theorem}

In particular in this case $\L^n(\ell^\infty_0,F)$ is pre-Riesz so has a Riesz completion for us to attempt to identify.

\begin{proposition}\label{ocod}
If $F$ is an almost Dedekind $\sigma$-complete Riesz space then\break $\L^n(\ell^\infty_0,F)$ is a majorizing and pervasive subspace of $\{T\in \L^n(\ell^\infty_0,F^\sigma):Te_n\in F\ \forall n\in\Natural\}$, and is hence order dense.
\begin{proof}
If $T\in\L^n(\ell^\infty_0,F^\sigma)$ with $Te_n\in F$ for all $n\in\Natural$, consider $T\1\in F^\sigma$. Pick $y\in F$ with $y\ge T\1$. As $F$ is almost Dedekind $\sigma$-complete, we may find a positive sequence $(z_n)$ in $F$ such that $\sum_{n=1}^\infty z_n \oto y-T\1$. Now define $S\in \L^r(\ell^\infty_0,F)$ with $S\1=y$ and $Se_n=Te_n+z_n$ for all $n\in\Natural$, using Proposition \ref{obexists}, and  note that
\[\sum_{k=1}^n Se_k=\sum_{k=1}^n Te_k+\sum_{k=1}^n z_k\oto T\1+(y-T\1)=S\1,\]
so that $S\in \L^n(\ell^\infty_0,F))$ by Theorem \ref{ocasDc}. It is simple to check that $S\ge T$, so that $\L^n(\ell^\infty_0,F)$ is majorizing.

To see that $\L^n(\ell^\infty_0)$ is pervasive, repeat the proof  of Theorem \ref{pervasive}, noting that the first case is impossible here and that in the second case the constructed operator $S\in \L^\diamond(\ell^\infty_0,F)=\L^n(\ell^\infty_0,F)$.

The deduction that $\L^n(\ell^\infty_0,F))$ is order dense is Proposition 2.8.5 of \cite{KvG}.
\end{proof}
\end{proposition}

\begin{proposition}\label{ocbandproj}If $F$ is an almost Dedekind $\sigma$-complete space then
\begin{enumerate}
\item $\{T\in\L^n(\ell^\infty_0,F^\sigma):Te_n\in F\ \forall n\in\Natural\}$ is a projection band in $\{T\in\L^r(\ell^\infty_0,F^\sigma):Te_n\in F\ \forall n\in\Natural\}$.
\item The band projection $P$ onto $\{T\in\L^n(\ell^\infty_0,F^\sigma)$ maps the operator $T$ to an operator $P(T)$ where $P(T)e_n=Te_n$ for all $n\in\Natural$ and $\sum_{k=1}^n Te_n\oto P(T)\1$.
\item The complementary band, $\{T\in\L^n(\ell^\infty_0,F^\sigma):Te_n\in F\ \forall n\in\Natural\}^d$, consists of all operators of the form $x=(x_n)\mapsto \lim_{n\to\infty}x_n\times y$ for $y\in F^\sigma$, and hence
\item $\{T\in \L^n(\ell^\infty_0,F^\sigma):Te_n\in F\ \forall n\in\Natural\}^d$ is linearly order isomorphic to $F^\sigma$.
\end{enumerate}
\begin{proof}
Recall that $F^\sigma$ is Dedekind $\sigma$-complete, so Theorem \ref{ocband} shows that\break $\L^n(\ell^\infty_0,F^\sigma)$ is a projection band in $\L^r(\ell^\infty_0,F^\sigma)$, and that the band projection has the form claimed. The restriction to $\{T\in \L^r(\ell^\infty_0,F^\sigma):Te_n \in F\ \forall n\in\Natural\}$ clearly has all the properties claimed.
\end{proof}
\end{proposition}

\begin{theorem}\label{ocRcompl}
If $F$ is an almost Dedekind $\sigma$-complete Riesz space then the Riesz completion of $\L^n(\ell^\infty_0,F)$, $\L^n(\ell^\infty_0,F)^\rho$ may be identified with the Riesz space of operators $\{T\in\L^n(\ell^\infty_0,F^\sigma):Te_n\in F\ \forall n\in\Natural\}$.
\begin{proof}

As $F^\sigma$  is Dedekind $\sigma$-complete, $\L^r(\ell^\infty_0,F^\sigma)$ is a Riesz space. By Theorem \ref{ocifDsc} $\L^n(\ell^\infty_0,F^\sigma)$ is a Riesz subspace of $\L^r(\ell^\infty_0,F^\sigma)$ and it is clear that $\{T\in \L^n(\ell^\infty_0,F^\sigma):Te_n\in F\ \forall n\in\Natural\}$ is a Riesz subspace  of that. The Proposition \ref{ocod} gives us the required order density. Finally, we need to show that if $T\in \L^n(\ell^\infty_0,F^\sigma)$ with $Te_n\in F$ for all $n\in\Natural$, then $T$ lies in the Riesz subspace generate by $\L^n(f,F)$. It clearly suffices to prove this when $T\ge 0$, in which case each $Te_n\ge 0$. Define $U\in\L^r(c,F)$ with
\[U\1=0\text{ and }Ue_n=\begin{cases}Te_n&\text{ if $n$ is odd}\\
-Te_{n-1}&\text{ if $n$ is even.}\\
\end{cases}\]
Proposition \ref{obexists} tells us, from the regularity of $T$, that $\big(\sum_{k=1}^n Te_k\big)_{n=1}^\infty$ is bounded above in $F^\sigma$ and hence in $F$. This is enough to give the existence of $U\in\L^r(c,F)$ and also that $Te_n\oto 0$. As $\sum_{k=1}^{2n} Ue_k=0$ and $\sum_{k=1}^{2n+1} Ue_k=Te_{2n+1}$ for each $n\in\Natural$, we have $\sum_{k=1}^n Ue_k\oto 0=U\1$, so that $U\in\L^n(c,F)$ by Theorem \ref{ocasDc}. By now familiar methods, \[U^+e_n=(Ue_n)^+=\begin{cases}Te_n&\text{ if $n$ is odd}\\
0&\text{ if $n$ is even}\\
\end{cases}\]
and $U^+\1=\sum_{k=1}^\infty Te_{2k-1}$. Similarly we may find $V\in\L^n(c,F)$ such that
\[V^+e_n=\begin{cases}0&\text{ if $n$ is odd}\\
Te_n&\text{ if $n$ is even}\\
\end{cases}\]
and $V^+\1=\sum_{k=1}^\infty Te_{2k}$. As $U^++V^+=T$, this completes the proof.
\end{proof}
\end{theorem}

\begin{corollary}\label{ascbandinRc}If $F$ is an almost Dedekind $\sigma$-complete Riesz space then\break $\L^n(\ell^\infty_0,F)^\rho$ is a projection band in $\L^r(\ell^\infty_0,F)^\rho$.
\end{corollary}

We refer the reader to \S4.2.1 of \cite{KvG} for the definition of a band in a pre-Riesz space. Corollary 4.2.7 of \cite{KvG} gives us:

\begin{corollary}\label{ascband}If $F$ is an almost Dedekind $\sigma$-complete Riesz space then $\L^n(\ell^\infty_0,F)$ is a band in $\L^r(\ell^\infty_0,F))$.
\end{corollary}

However this band need not be a projection band see \cite{KM}. In Theorem 14 of \cite{KM} it is shown that a linear operator $P$ on a pre-Riesz space $E$ is a band projection if and only if it is an order projecion, i.e. $0\le P\le I_E$, where $I_E$ is the identity on $E$.

\begin{proexample}\label{bandnotproj}$\L^n(\ell^\infty_0)$ is a band but not a projection band in $\L^r(\ell^\infty_0)$.
\begin{proof}
As $\ell^\infty_0$ is almost Dedekind $\sigma$-complete, Corollary \ref{ascband} tells us that $\L^n(\ell^\infty_0)$ is a band in $\L^r(\ell^\infty_0)$.
Using Propositions \ref{obexists} and \ref{oborder}, there are operators $T$ and $T_m$ in $\L^r(\ell^\infty_0)_+$ for $m\in\Natural$, with (i) $T\1=\1$ and $Te_n=e_{2n}$ for $n\in\Natural$ and (ii) $T_m\1=\sum_{k=1}^m e_{2k}$ and $T_m e_n=e_{2n}$ if $n\le m$ else 0. By Theorem \ref{ocasDc} each $T_m\in\L^n(\ell^\infty_0)$. It is routine to check that $0\le T_m\le T$ for all $m\in\Natural$. Suppose that $P$ were the band projection of $\L^r(\ell^\infty_0)$ onto $\L^n(\ell^\infty_0)$.
 Then for each $n\in\Natural$, $0\le P(T)e_n\le Te_n=e_{2n}$ so that $P(T)e_n=\alpha_{2n}e_{2n}$ for $\alpha_{2n}\in [0,1]$. As $\sum_{k=1}^n P(T)e_n\uparrow P(T)\1$, and odd entries in each $\sum_{k=1}^n P(T)e_n$ are all zero, the eventual constant value of $P(T)\1$ is zero. Hence $(\alpha_{2n})$ is eventually zero. On the other hand $\alpha_{2n}e_{2n}=P(T)e_n\ge P(T_n)e_n=T_n e_n=e_{2n}$ for all $n\in\Natural$, which is a contradiction.
\end{proof}
\end{proexample}

\begin{problem}Either remove the hypothesis that $F$ is almost Dedekind $\sigma$-complete from recent results, or give an example to show that it cannot be removed.
\end{problem}
\section{Attempts to extend the domain to be $c$}

An order bounded operator $T$ from $\ell^\infty_0$ into a uniformly complete Riesz space $F$ certainly maps the whole of $\ell^\infty_0$ into a principal ideal in $F$, say $F_u$. $F_u$ is complete for the order unit norm on $F_u$ induced by $u$, $\|\cdot\|_u$. If we give $\ell^\infty_0$ the usual supremum norm, $\|\cdot\|_\infty$ then $T$ is certainly bounded for these two norms, so is continuous.Now $T$ may be extended by continuity to the whole of $c$. It follows that many of our results, in particular \ref{regob}, \ref{obexists}, \ref{oborder}, \ref{pospart}, \ref{Rcompl}, \ref{posoc}, \ref{ocifpospart}, \ref{ocifDsc}, \ref{ocband}, \ref{ocasDc}, \ref{ocpg}, \ref{ocRcompl}, \ref{ascbandinRc}, \ref{ascband} and \ref{bandnotproj}  remain valid if, in their statements, $\ell^\infty_0$ is replaced by $c$ and $F$ is specified to be uniformly complete, with very little (if any) extra detail being needed for their proofs. The only additional argument needed is in the proof of Theorem 2.1 where the earliest reference that we could find for the statement that $\L^r(c,F)$ is a Riesz space if and only if $F$ is Dedekind $\sigma$-complete is Lemma 4.14 of \cite{vR}.

There are notable gaps in our knowledge if we do not assume that $F$ is uniformly complete. In particular we ask:

\begin{question}
Is there a description of the Riesz completion of $\L^r(c,F)$, when $F$ is only assumed to be an Archimedean Riesz space, similar to that in Theorem \ref{Rcompl}?
\end{question}

Theorem 3.3 of \cite{AW} shows that $\L^r(c,\ell^\infty_0)$ is a Riesz space, so we cannot expect as precise a result as our Theorem \ref{Rcompl} but it might not be too much to expect a description of the form $\{T\in\L^r(c,G):Te_n\in F\ \forall n\in\Natural\}$ where $G$ is a Riesz subspace of $F^\sigma$ containing $F$.

Even with such a result there are obstacles to describing the Riesz completion of $\L^n(c,F)$. Even with the assumption of being almost Dedekind $\sigma$-complete, we don't yet have a description of the order continuous operators. It may be worth recording that there actually is a description similar to our previous results. We  know of no reference for the following (almost) obvious lemma.

\begin{lemma}\label{oconvcompl}Let $E$ be an Archimedean Riesz space with a strong order unit $u$ and $\overline{E}$ be its completion for the order unit norm, $\|\cdot\|_u$, induced by $u$. A net $(a_\gamma)_{\gamma\in\Gamma}$ is order convergent to 0 in $E$ if and only if it is order convergent to 0 in $\overline{E}$.
\begin{proof}
If $a_\gamma\oto 0$ in $E$ then certainly $a_\gamma\oto 0$ in $\overline{E}$. Suppose that $a_\gamma\oto 0$ in $\overline{E}$, so there is a net $(b_\lambda)_{\lambda\in\Lambda}$ in $\overline{E}$ decreasing to 0 and such that for all $\lambda\in\Lambda$ there is $\gamma_0$ such that $|a_\gamma|\le b_\lambda$ for all $\gamma\ge \gamma_0$. We need to replace the net $(b_\lambda)$ in $\overline{E}$ by a net in $E$. For each $n\in\Natural$ and $\lambda\in\Lambda$, let $b_{\lambda,n}\in E$ with $b_\lambda\le b_{\lambda,n}\le b_\lambda+u/n$. Certainly the net $(b_{\lambda,n})_{\Gamma\times\Natural}$ lies in $E$ and has infimum 0, but may not be decreasing, so we take our net to be indexed by $n$-tuples of arbitrary length from $\Lambda\times \Natural$, with
\[c_{((\lambda_1,n_1), (\lambda_2,n_2),\dots,(\lambda_m,n_m))}=\bigwedge_{k=1}^m b_{\lambda_k, n_k},\]
which now does decrease to 0. Given an index $((\lambda_1,n_1),\dots,(\lambda_m,n_m))$ there are $\gamma_k$ such that $|a_\gamma|\le b_{\lambda_k}$ if $\gamma\ge \gamma_k$. Now take $\gamma_0\ge \gamma_1,\gamma_2,\dots, \gamma_m$. If $\gamma\ge \gamma_0$ then for each $k=1,2,\dots,m$, $|a_\gamma|\le b_{\lambda_k}\le b_{\lambda_k,n_k}$, so that
\[|a_\gamma|\le \bigwedge_{k=1}^m b_{\lambda_k,n_k}=c_{(\lambda_1,n_1),\dots,(\lambda_m,n_m))}.\]
\end{proof}
\end{lemma}

\begin{theorem}
Let $F$ be an almost Dedekind $\sigma$-complete Riesz space and $T\in \L^b(c,F)$, then the following are equivalent:
\begin{enumerate}
\item $T$ is order continuous.
\item $T$ is $\sigma$-order continuous.
\item $T\1-\sum_{k=1}^n Te_k\oto 0$.
\end{enumerate}
\begin{proof}
It is only (3) implies (1) that we need to prove. As $T$ is order bounded $T(c)$ is contained in a principal ideal $F_u$ in $F$. Consider $T$ as an operator from $c$  into the $\|\cdot\|_u$ completion of $F_u$. It is simple to see that this completion is almost Dedekind $\sigma$-complete if $F$ is. If $T\1-\sum_{k=1}^n Te_k\oto 0$ in $E$ then certainly $T\1-\sum_{k=1}^n Te_k\oto 0$ in $\overline{F_u}$. As $\overline{F_u}$ is uniformly complete, $T\in\L^n(c,\overline{F_u})$. If $(x_\gamma)$ is a net in $c$ with $x_\gamma\oto 0$ then $Tx_\gamma\oto 0$ in $\overline{F_u}$ as $T\in\L^n(c,\overline{F_u})$. By Lemma \ref{oconvcompl} $Tx_\gamma\oto 0$ in $F_u$, and hence in $F$, so that $T\in\L^n(c, F)$.
\end{proof}
\end{theorem}

Unfortunately, even with this result we do not know that $\L^n(c,F)$ is positively generated so it may  not be pre-Riesz.

\begin{question}
Is $\L^n(c,F)$ positively generated for all Archimedean Riesz spaces~$F$? Recall that order boundedness is part of our definition of an operator being in  $\L^n(c,F)$. If it makes a difference, an answer for the space of regular order continuous operators would be of interest.
\end{question}

One important case of  $F$ being uniformly complete is when it is a Banach lattice. In this case we have the regular norm on $\L^r(c,F)$ and also the regular norm on our concrete realization of its Riesz completion. An obvious question to ask is now the two are related. First we need to put a norm on $F^\sigma$. This is easy, though by no means unique. Indeed, Lemma 3.2.1 of \cite{KvG} tells us that defining $\|y\|=\inf\{\|x\|:x\in F,\ x\ge|y|\}$, the \emph{maximal} extension, extends the Riesz norm on $F$ to a Riesz norm on $F^\sigma$. There is, of course, no guarantee that $F^\sigma$ will be complete for this norm.

\begin{theorem}
Let $F$ be a Banach lattice. On $\L^r(c,F)^\rho=\{T\in\L^r(c,F^\sigma):T(c_0)\subset F\}$, the regular norm is a Riesz norm that extends the regular norm on $\L^r(c,F)$. Hence the regular norm on $\L^r(c,F)$ is a pre-Riesz norm (\cite{KvG}, Definition 3.2.24.)
\begin{proof}
We need to show that if $T\in \L^r(c,F)$ then its regular norm, $\|T\|_1$, is the same as the regular norm, $\|T\|_2$, when we consider it as an element of $\L^r(c,F)^\rho$. As an  element of $\L^r(c,F)^\rho$ it has a modulus and  $\|T\|_2$ is just the operator norm of $|T|$. As $c$ has a strong order unit, $\1$, this operator norm is precisely $\big\||T|\1\big\|$. By definition of the extension of the norm on $F$ to $F^\sigma$, $\big\||T|\1\big\|=\inf\{\|x\|:x\in F, x\ge |T|\1\}$. For each such $x$, there is $T_x\in\L^r(c,F)$ with $T_xe_n=|T|e_n$ for all $n\in\Natural$ and $T_x\1=x$. It is routine to check that $T_x\ge |T|\ge \pm T$, so that $\|T\|_1$ is at most the operator norm of $T_x$, which is precisely $\|x\|$. Hence $\|T\|_1\le \|T\|_2$. As the reverse inequality is clear, we have $\|T\|_1=\|T\|_2$.

The final statement follows from
Theorem 3.2.31 of \cite{KvG}.
\end{proof}
\end{theorem}

If $F$ is, in addition, almost Dedekind $\sigma$-complete then we can say more. We simply need the fact that $F^\sigma$ is a Banach lattice in this case. We know no convenient reference for this explicit result. It seems to be a direct consequence of \cite{VL} but the author does not have access to it. In view of this, and the difficulty that readers may have in accessing it, we give the short proof here.

\begin{theorem}\label{menc}If $F$ is an almost Dedekind $\sigma$-complete Banach lattice then $F^\sigma$ under the maximal extension norm is norm complete.
\begin{proof}
In order to prove that $F^\sigma$ is norm complete, it suffices to show that every absolutely convergent series in $F^\sigma$ is convergent. If $\sum_{n=1}^\infty \|x_n\|<\infty$ then also $\sum_{n=1}^\infty \|x_n^\pm\|<\infty$, so it suffices to consider an absolutely convergent series of positive terms, $\sum_{n=1}^\infty x_n$. By definition of the maximal norm, we may find $y_n\in F$ with $0\le x_n\le y_n$ and $\|y_n\|\le 2\|x_n\|$. Thus $\sum_{n=1}^\infty y_n$ is absolutely convergent, and hence convergent in $F$. Thus $\|\sum_{n=m}^\infty y_n\|\to 0$ as $m\to\infty$. We know that $F^\sigma$ is Dedekind $\sigma$-complete, so $\sum_{n=1}^\infty x_n\oto x\in F^\sigma$. As the order convergent sum $\sum_{n=m}^\infty x_n$ is smaller than the norm convergent sum $\sum_{n=m}^\infty y_n$, we have $\|\sum_{n=m}^\infty x_n\|\to 0$ as $m\to\infty$. I.e.
$\|\sum_{n=1}^{m-1} x_n-x\|\to 0$ as $m\to\infty$.
\end{proof}
\end{theorem}

\begin{corollary}
If $F$ is an almost Dedekind $\sigma$-complete Banach lattice then the regular norm on $\L^r(c,F)$ is a pre-Riesz norm that extends to a Banach lattice norm on $\L^r(c,F)^\rho$.
\end{corollary}

In particular the last conclusion holds if $F$ is a separable Banach lattice.

\begin{question}Is the maximal extension norm on $F^\sigma$ norm complete for every Banach lattice $F$?
\end{question}

\section{Miscellaneous observations}

It is worth recording one special example. Recall the classical description of infinite matrices which act on $c$.

\begin{theorem}[Kojima-Schur]\label{KS} The infinite matrix $A=(a_{mn})$ transforms $c$ into itself via matrix multiplication, $x=(x_n)\mapsto (\sum_{n=1}^\infty a_{mn} x_n)_{m=1}^\infty$ if and only if  the following three conditions hold:
\begin{enumerate}
\item[(KS1)] $\{\sum_{n=1}^\infty |a_{mn}|:m\in\Natural\}$ is bounded.
\item[(KS2)] $\lim_{m\to\infty} a_{mn}$ exists for all $n\in\Natural$.
\item[(KS3)] $\lim_{m\to\infty}\sum_{n=1}^\infty a_{mn}$ exists.
\end{enumerate}
\end{theorem}

A proof may be found in \cite{H}, \S3.2, Theorem 1. By considering the image of standard basic vectors it is clear that the induced operator is positive if and only if each $a_{mn}\ge 0$, i.e. $A$ is a positive matrix.
Perusal of the proof in \cite{H} makes it clear that condition (KS1) states precisely that $A$ maps $c$ into $\ell^\infty$. It is also evident that (KS2) is equivalent to $A$ mapping standard basic vectors into $c$, from which it is a small step to concluding that $A$ maps $c_0$ into $c$. Condition (KS3) is what is needed for $\1$ to be mapped into $c$. It should be apparent that the operators induced by these matrices are precisely the operators in $\L^\diamond(c)$ and hence, as $c$ is almost Dedekind $\sigma$-complete, in $\L^n(c)$. Therefore, as $c^\sigma=\ell_\infty$, we have:

\begin{theorem}The Riesz completion of $\L^n(c)$ may be identified with the infinite matrices $A$ satisfying (KS1) and (KS2) of Theorem \ref{KS}.
\end{theorem}

Finally, let us point out that there will be versions of most of our results for analogues of $\ell^\infty_0$ and $c$ involving larger numbers of atoms. It is not clear to the author that anything fundamentally new will be learnt from them.

\end{document}